\documentclass[a4paper,11pt,twoside,reqno]{amsart}

\usepackage[utf8]{inputenc}
\usepackage[plainpages=false,pdfpagelabels=true]{hyperref}
\usepackage{amssymb,amsthm}
\usepackage[margin=1in]{geometry}
\usepackage{slashed}
\usepackage{comment}
\usepackage{caption}

\usepackage[toc,page]{appendix}

\newtheorem{Satz}{Theorem}[section]
\newtheorem{Prop}[Satz]{Proposition}
\newtheorem{Lem}[Satz]{Lemma}

\newtheorem{Thm}[Satz]{Theorem}
\newtheorem{Cor}[Satz]{Corollary}
\theoremstyle{definition}
\newtheorem{Dfn}[Satz]{Definition}
\newtheorem{Bem}[Satz]{Remark}

\newcommand{\tr}{\operatorname{Tr}}
\newcommand{\sff}{\mathrm{I\!I}}
\newcommand{\hess}{\operatorname{Hess}}

\newcommand{\vol}{{\operatorname{vol}}}
\newcommand{\dv}{\text{ }dv}

\newcommand{\norm}[1]{\left \lVert #1 \right \rVert}
\allowdisplaybreaks[1]

\renewcommand{\epsilon}{\varepsilon}

\newcommand{\R}{\ensuremath{\mathbb{R}}}
\newcommand{\N}{\ensuremath{\mathbb{N}}}

\newcommand{\s}{\ensuremath{\mathbb{S}}}
\numberwithin{equation}{section}

\usepackage{color}

\newtheorem{innercustomgeneric}{\customgenericname}
\providecommand{\customgenericname}{}
\newcommand{\newcustomtheorem}[2]{%
  \newenvironment{#1}[1]
  {%
   \renewcommand\customgenericname{#2}%
   \renewcommand\theinnercustomgeneric{##1}%
   \innercustomgeneric
  }
  {\endinnercustomgeneric}
}

\newcustomtheorem{customthm}{Theorem}
\newcustomtheorem{customcor}{Corollary}

\title{On the stability of the generalized equator map}
\author{Volker Branding}
\date{\today}
\address{University of Vienna, Faculty of Mathematics\\
Oskar-Morgenstern-Platz 1, 1090 Vienna, Austria\\}
\email{volker.branding@univie.ac.at}

\author{Anna Siffert}
\address{Universität M\"unster, Mathematisches Institut\\
Einsteinstr. 62\\
48149 M\" unster\\
Germany}
\email{asiffert@uni-muenster.de}

\subjclass[2010]{58E20; 53C43}
\keywords{generalized equator map; extrinsic polyharmonic map; stability}
\thanks{The first author gratefully acknowledges the support of the Austrian Science Fund (FWF) through the project "Geometric Analysis of Biwave Maps" (DOI: 10.55776/P34853).}

\begin{document}
\begin{abstract}
The energy, the $p$-energy ($p\in\mathbb{R}$ with $p\geq 2$) and the extrinsic $k$-energy ($k\in\mathbb{N}$) for maps between Riemannian manifolds are central objects in the geometric calculus of variations. The equator map from the unit ball to the Euclidean sphere provides an explicit critical point of all aforementioned energy functionals. During the last four decades many researchers studied the stability of this particular map when considered as a critical point of one of these energy functionals, see e.g. \cite{MR4436204}, \cite{MR705882}. 

Recently, Nakauchi \cite{MR4593065} introduced a generalized radial projection map
and proved that this map is
both a critical point of the energy and a critical point of the $p$-energy. This generalized radial projection map gives rise to a generalized equator map which is also both a critical point of the energy and a critical point of the $p$-energy.

In this manuscript we first of all show that the generalized equator map 
is also a critical point of the extrinsic $k$-energy. Then, the main focus is a detailed
stability analysis of this map, considered as a critical point of both the extrinsic $k$-energy and the \(p\)-energy.
We thus establish a number of interesting generalizations of the classical (in)stability results of Jäger and Kaul \cite{MR705882}.
\end{abstract} 

\maketitle

\section{Introduction}
In recent years the study of higher order
energy functionals has become an increasingly important research area within geometric analysis, see e.g. the recent manuscripts \cite{MR4106647, MR4817500}.
There exist several different definitions of energy functionals, some being natural from the geometric point of view and others are tailored for the use of analytic methods. 

\smallskip
Concerning the energy functionals that are natural from a geometric point of view we mention the well-studied energy functional
\begin{align}
\label{harmonic}
E(\phi)=\frac{1}{2}\int_M|d\phi|^2\dv,    
\end{align}
where \(\phi\colon M\to N\) is a smooth map from a compact Riemannian manifold
\((M,g)\) to a Riemannian manifold \((N,h)\). 
The critical points of \eqref{harmonic} are characterized by the vanishing of the so-called \emph{tension field}

\begin{align*}
    0=\tau(\phi):=\tr\bar\nabla d\phi,\qquad \tau(\phi)\in\Gamma(\phi^\ast TN),
\end{align*}
where \(\bar\nabla\) represents the connection on the pull-back bundle \(\phi^\ast TN\). Solutions of \(\tau(\phi)=0\) are called \emph{harmonic maps}.
For more details on harmonic maps we refer to the book \cite{MR2044031}.

\smallskip

There is also a second, extrinsic, point of view on harmonic maps. Here, one uses the Nash embedding theorem to isometrically embedded \(N\) into some $\mathbb{R}^q$ for $q\in\mathbb{N}$ sufficiently large.
In this case we consider maps \(u\colon M\to\R^q\) and the extrinsic version of 
\eqref{harmonic} is simply
\begin{align}
\label{eq:e-energy}
E^{ext}_1(u)=\frac{1}{2}\int_M|du|^2dv
\end{align}
with critical points given by solutions of the partial differential equation
\begin{align*}
\Delta u+\sff(u)(du,du)=0,    
\end{align*}
where \(\sff(u)\colon M\to\R^q\) is the second fundamental form of the embedding.

The energy functional (\ref{eq:e-energy}) is coercive but it depends on the embedding of the target manifold into Euclidean space. 
This makes the extrinsic version \eqref{eq:e-energy} of the energy functional suitable for analytic considerations but not favorable from a geometric point of view.

In this manuscript we will study specific maps from Euclidean balls to spheres exclusively, thus we introduce the extrinsic higher order energy functionals just for spherical targets: 
Let $u:M\rightarrow \s^{n}$ be a map and let $\iota:\s^{n}\rightarrow\mathbb{R}^{n+1}$ be the canonical embedding. 
We will deal with maps in the Sobolev space
\begin{align*}
W^{k,2}(M,\s^{n})=\{ u\in W^{k,2}(M,\mathbb{R}^{n+1})\,\colon\, u(x)=(u_1(x),\dots, u_{n+1}(x))\in\s^{n}\,\mbox{almost everywhere}\},   
\end{align*}
where $k\in\mathbb{N}$ and $(u_1,\dots,u_{n+1})$
denote the component functions of $\iota\circ u$.
The \textit{extrinsic $k$-energy functional} for $u\in W^{k,2}(M,\s^{n})$ is defined by
\begin{align}
\label{poly-energies}
E_k^{ext}(u)&=\frac{1}{2}\int_{M}|\Delta^su|^2\dv,\qquad\mbox{when}\,\, k=2s, s\in\mathbb{N};\\
\nonumber E_k^{ext}(u)&=\frac{1}{2}\int_M|\nabla\Delta^su|^2\dv,\qquad\mbox{when}\,\,  k=2s+1,  s\in\mathbb{N}_0.
\end{align}
Critical points of the extrinsic $k$-energy functional are called \textit{extrinsic $k$-harmonic maps} or \textit{extrinsic polyharmonic maps} of order \(k\).
An important class of critical points of the extrinsic $k$-energy functional is constituted by the so-called \textit{minimizers (of the extrinsic $k$-energy functional)}. More precisely, a map $u\in W^{k,2}(B^m,\s^{n})$ is called minimizer if $E_k^{ext}(u)\leq E_k^{ext}(v)$ for all $v\in W^{k,2}(B^m,\s^{n})$ with $u-v\in W_0^{k,2}(B^m,\mathbb{R}^{n+1})$.
Furthermore, an extrinsic $k$-harmonic map is called \textit{stable} if the second variation of \eqref{poly-energies}
is positive, i.e.
\begin{align*}
    \frac{d^2}{dt^2}E_k^{ext}(u_t)\lvert_{t=0}\geq 0
\end{align*}
for all variations $u_t$ with $u_t-u\in W_0^{k,2}(B^m,\mathbb{R}^{n+1})$. Otherwise the map $u$ is called \textit{unstable}. Note that we do not distinguish between stable and weakly stable critical points in this manuscript.

For the sake of completeness we mention that
one can also define extrinsic polyharmonic maps as critical points of
\begin{align*}
\tilde E^{ext}_k(u)=\frac{1}{2}\int_M|\nabla^ku|^2\dv,\qquad k\in\N,
\end{align*}
see for example \cite{MR2267035,MR2520907} for some analytic results
on such maps.
\medskip

One specific critical point of the extrinsic energy functionals which has been studied intensively in the literature is the
equator map, i.e. the map
\begin{align}
\label{equator}
  u_\star^{(1)}\colon B^m&\to\s^m\subset\mathbb{R}^m\times\mathbb{R},\\ 
  \notag x&\mapsto\big(\frac{x}{r},0\big),
\end{align}
where $r=r(x):=\norm{x}$ denotes the Euclidean distance to the origin.
The equator map is contained in $W^{k,2}(B^m,\s^m)$ if and only if $m\geq 2k+1$.
Furthermore, if $m\geq 2k+1$, it is an extrinsic $k$-harmonic map, see 
\cite[Proposition 1.1]{MR4436204}. 

Concerning the stability of \eqref{equator} we want to mention the seminal
result of Jäger and Kaul \cite{MR705882} who showed that the equator map is minimizing if $m\geq 7$ and unstable if $3\leq m<7$ when considered as a harmonic map.
Later, in \cite{MR2322746} Hong and Thompson studied the stability of (\ref{equator}) as an extrinsic biharmonic map, i.e. $k=2$, and showed that the equator map
is minimizing if $m\geq 10$ and unstable if $5\leq m<10$.
Recently, Fardoun, Montaldo and Ratto extended this stability analysis to the case of the equator map as a critical point of the extrinsic $k$-energy for $k\geq 3$
in \cite{MR4436204}. In order to achieve their results the authors conduct a sophisticated analysis of the sign of
a polynomial $P_k(m)$ (see \cite{MR4436204} or Section\,\ref{sec:prelim} for the definition of $P_k(m)$). They prove that for $m\geq 2k+1$ the equator map is minimizing if $P_k(m)\geq 0$ and unstable if $P_k(m)<0$, see \cite[Theorem 1.4]{MR4436204}.  

\medskip

In \cite{MR4371934,MR4593065} a \lq generalized radial projection\rq\, has been introduced.
Namely, in \cite[Main Theorem, p.1]{MR4593065} Nakauchi showed that for any $\ell,m\in\mathbb{N}$ with $\ell\leq m$ there exists a harmonic map
\begin{align}
\label{nak-maps}
u^{(\ell)}\colon\mathbb{R}^{m}\setminus\{0\}&\to\s^{m^{\ell}-1},\\
\notag x=(x_1,\dots,x_m)&\mapsto u^{(\ell)}(x)=(u^{(\ell)}_{i_1\dots i_{\ell}}(x))_{1\leq i_1,\dots i_{\ell}\leq m}.
\end{align}
We set \(y_i=\frac{x_i}{r}\) and have the recursive definition
\begin{align*}
u^{(1)}_{i_1}(x)=&y_{i_1},\\
u^{(\ell)}_{i_1\ldots i_\ell}(x)=&C_{\ell,m}\big(y_{i_1}u^{(\ell-1)}_{i_1\ldots i_{\ell-1}}(x)-\frac{1}{\ell+m-3}r\frac{\partial}{\partial x_{i_\ell}}u^{(\ell-1)}_{i_1\ldots i_{\ell-1}}(x)\big),
\end{align*}
where 
\begin{align*}
 C_{\ell,m}=\sqrt{\frac{\ell+m-3}{2\ell+m-4}},   
\end{align*}
leading to the generalized equator map 

\begin{align}
\label{equator-generalized}
  u_\star^{(l)}\colon B^m&\to\s^{m^{\ell}}\subset\mathbb{R}^{m^{\ell}}\times\mathbb{R},\\ 
  \notag x&\mapsto\big(u^{(l)},0\big),
\end{align}    
where \(u^{(l)}\) is defined in \eqref{nak-maps}.

For $\ell=1$ the equator map (\ref{equator}) is recovered.
Further, note that we consider $u_\star^{(\ell)}$ as one-parameter family of maps and use the singular to refer to it.

\medskip

It is natural to study the stability of the generalized equator map \eqref{equator-generalized}, which is the content of this manuscript.
In the following, we will present the main results of this article.
Our first result concerning the map \(u_\star^{(\ell)}\) is the following:
\begin{customthm}{\ref{thm:sobolev}}
Let $\ell,m\in\mathbb{N}$ with $\ell\leq m$.
The generalized equator map \(u_\star^{(\ell)}\colon B^m\to\s^{{m^\ell}}\) defined
in \eqref{equator-generalized} belongs to \(W^{k,2}(B^m,\s^{{m^\ell}})\) if and only if \(m\geq 2k+1\). In particular, \(u_\star^{(\ell)}\) is an extrinsic $k$-harmonic map.
\end{customthm}

Afterwards, we study the stability of the map \eqref{equator-generalized} 
considered as extrinsic $k$-harmonic map, $k\in\mathbb{N}$,
and obtain the following result, which widely generalizes \cite[Theorem 1.4]{MR4436204}:

\begin{customthm}{\ref{thm:stability}}
Assume \(m\geq 2k+1\), where $m,k\in\mathbb{N}$. 
There exists a polynomial \(Q_k^\ell(m)\) such that the generalized equator map \(u_\star^{(\ell)}\colon B^m\to\s^{m^{\ell}}\) is
\begin{enumerate}
    \item an energy minimizing extrinsic \(k\)-harmonic map if \(Q_k^\ell(m)\geq 0\),
    \item an unstable extrinsic \(k\)-harmonic map if \(Q_k^\ell(m)<0\).
\end{enumerate}
\end{customthm}

Note that for $\ell=1$ Theorem\,\ref{thm:stability} reduces to \cite[Theorem 1.4]{MR4436204}. Further, Theorem\,\ref{thm:stability} implies the following Corollary:

\begin{customcor}{\ref{cor:harmonic}}
The generalized equator map \(u_\star^{(\ell)}\colon B^m\to\s^{{m^\ell}}\) defined
in \eqref{equator-generalized} is 
\begin{enumerate}
    \item an energy minimizing harmonic map if
    \begin{align*}
m\geq 2(l+1)+2\sqrt{2}\ell;    
\end{align*}
\item an unstable harmonic map if 
\begin{align*}
    3\leq m <2(l+1)+2\sqrt{2}\ell.
\end{align*}
\end{enumerate}
\end{customcor}

Note that for $\ell=1$ Corollary \ref{cor:harmonic} recovers the famous (in)stability
result of Jäger and Kaul \cite{MR705882}.

\smallskip

A careful study of the sign of \(Q_k^\ell(m)\) in dependence of $k,\ell$ and $m$, yields the following Theorem:

\begin{customthm}{\ref{con}}
For any $k\geq 2$ and $\ell\in\mathbb{N}$ we have $Q_k^{\ell}(m)>0$ for any $m\geq 5+3(k-2)+5\ell$. If this inequality is satisfied, the generalized equator map $u_\star^{(\ell)}:B^m\rightarrow\s^{m^{\ell}}$ is an energy minimizing extrinsic $k$-harmonic map.
\end{customthm}

We also discuss the case \(\ell=m\) in the stability analysis:

\begin{customthm}{\ref{thm:lm}}
Let \(Q_k^\ell(m)\) be the polynomial defined in \eqref{dfn:Q}. Assume that \(\ell=m\geq 2k+1\). Then, the generalized equator map \(u_\star^{(m)}\colon B^m\to\s^{m^{m}}\) is unstable. 
\end{customthm}

Another possible higher order generalization of the well-studied energy functional \eqref{harmonic}
for maps between Riemannian manifold is given by the \(p\)-energy functional
which is defined by
\begin{align}
\label{p-energy}
E_p(\phi):=\frac{1}{p}\int_M|d\phi|^p\dv,
\end{align}
where \(p\in\R\) with \(p\geq 2\).
The critical points of \eqref{p-energy} are given by the solutions of the equation
\begin{align*}
0=\tr_g\bar\nabla(|d\phi|^{p-2}d\phi)
\end{align*}
and are called \emph{p-harmonic maps}.
In the case of maps to the sphere \(u\colon M\to\s^n\) this equation simplifies to
\begin{align}
\label{p-tension-field-sphere}
\tr_g\nabla(|\nabla u|^{p-2}du)+|\nabla u|^pu=0.    
\end{align}

One motivation for studying \eqref{p-energy} is the following:
It is well-known that the standard harmonic map energy functional \eqref{harmonic}
is invariant under conformal transformations if the domain is two-dimensional, which no longer holds true if the dimension of the domain is larger than two.
However, the $p$-energy is conformally invariant in certain higher dimensions,
more precisely, \eqref{p-energy} is conformally invariant in the case that \(p=\dim M\).
The terms stable, unstable and minimizing for $p$-harmonic maps are defined analogously to the $k$-extrinsic setting.

For more details on \(p\)-harmonic maps between Riemannian manifolds
we refer to \cite{MR4579395} and references therein.

\smallskip

The map \(u^{(\ell)}\) defined in \eqref{equator-generalized} is also a critical point of \eqref{p-energy}, i.e. it solves \eqref{p-tension-field-sphere},
see the appendix of \cite{MR4757446} for the precise details.
Here, we present a complete analysis of the stability of \eqref{equator-generalized} when considered as a \(p\)-harmonic map extending the results from \cite{MR4757446} and \cite{MR4686854}.

\begin{customthm}{\ref{thm:p-harmonic}}
Assume that \(m>p\geq 2\).
\begin{enumerate}
    \item If
     \begin{align*}
     \frac{(m-p)^2}{4\ell(\ell+m-2)} < 1
     \end{align*}
     then the generalized equator map \(u_\star^{(\ell)}\colon B^m\to\s^{m^\ell}\) is an unstable \(p\)-harmonic map.
     \item If
     \begin{align*}
        \frac{(m-p)^2}{4\ell(\ell+m-2)}\geq  1
     \end{align*}
     then the generalized equator map \(u_\star^{(\ell)}\colon B^m\to\s^{m^\ell}\) is the unique energy minimizer 
     of the \(p\)-energy \eqref{p-energy}.
\end{enumerate}
\end{customthm}

For \(\ell=1\) we reobtain \cite[Theorem A]{MR1743430}, where it was shown that the equator map \(u_\star^{(1)}\) is a minimizer of the \(p\)-energy \eqref{p-energy} if \(m\geq 2+p+2\sqrt{p}\), see also \cite{MR1647881}.

\medskip

\textbf{Notation and Conventions:}
Throughout this article we will employ the following sign conventions: 
For the Riemann curvature tensor we use 
$$
R(X,Y)Z=[\nabla_X,\nabla_Y]Z-\nabla_{[X,Y]}Z,
$$ 
where \(X,Y,Z\) are vector fields.

For the rough Laplacian on the pull-back bundle $\phi^{\ast} TN$ we employ the analysts sign convention, i.e.
$$
\bar\Delta = \tr(\bar\nabla\bar\nabla-\bar\nabla_\nabla).
$$
In particular, this implies that the Laplace operator has a negative spectrum.

By $dv$ we represent the volume element of an arbitrary Riemannian manifold $(M,g)$.
In the specific case that $(M,g)$ is the Euclidean ball, the volume element is denoted by $dx$.

\medskip

\textbf{Organization:}
Section\,\ref{sec:prelim} contains preliminaries and establishes Theorem\,\ref{thm:sobolev}.
 We investigate the stability of the generalized equator map in Section\,\ref{sec-stability}, in particular Theorem\,\ref{thm:stability} and Corollary\,\ref{cor:harmonic} are proven in this section. The sign of $Q_k^{\ell}(m)$
 in dependence of $k,\ell,m\in\mathbb{N}$
 is studied in Section\,\ref{sec-sign} which in particular contains Theorem\,\ref{con} and Theorem\,\ref{thm:lm}. Finally, Section\,\ref{sec-stability-p}
investigates the stability of the generalized equator map
when considered as a critical point of the \(p\)-energy. This section in particular contains Theorem\,\ref{thm:p-harmonic}.

\section{Preliminaries}
\label{sec:prelim}
In this section we provide basic facts on extrinsic polyharmonic maps to spheres and the generalized equator map introduced by Nakauchi \cite{MR4593065}.

\subsection{Extrinsic polyharmonic maps to spheres}
We briefly recall the Euler-Lagrange equations for extrinsic polyharmonic maps to spheres.
The critical points $u\in W^{k,2}(B^m,\s^{q-1})$ 
of \eqref{poly-energies} are characterized by the equation
\begin{align}
\label{eq:extrinsic-polyharmonic}
\Delta^ku-\langle\Delta^ku,u\rangle u=0.
\end{align} 
Exploiting the fact that we are considering a spherical target,
this equation can also be expressed in the form
\begin{align*}
\Delta^ku+\bigg(
\Delta^{k-1}(|\nabla u|^2)+\sum_{j=0}^{k-2}\Delta^j(\langle\Delta^{k-1-j}u,\Delta u\rangle)
+2\sum_{j=0}^{k-2}\Delta^j(\nabla\Delta^{k-1-j}u,\Delta u\rangle)
\bigg)u=0,
\end{align*}
see \cite[Proposition 2.1]{MR4436204} for the precise details.

\subsection{Generalized equator map}

We will now give some more background on the map defined in \eqref{nak-maps}. Let $\ell,m\in\mathbb{N}$ such that $\ell\leq m$.
Nakauchi \cite{MR4593065} proved that the map $u^{(\ell)}\colon\mathbb{R}^{m}\setminus\{0\}\to\s^{m^\ell-1}$ has the following properties: 
\begin{enumerate}
\item $u^{(\ell)}$ satisfies the equation for harmonic maps to spheres
\begin{align*}
    \Delta u^{(\ell)}+\lvert\nabla u^{(\ell)}\rvert^2 u^{(\ell)}=0; 
\end{align*}
\item $u^{(\ell)}$ is a polynomial in $u_{i_1},\dots, u_{i_{\ell}}$ of degree $\ell$, where $u_{i_j}=\frac{x_{i_j}}{r}$;
\item $\lvert\nabla u^{(\ell)}\rvert^2=\frac{\ell(\ell+m-2)}{r^2}$.
\end{enumerate}
Furthermore, Nakauchi \cite[Proposition 1, (1)]{MR4593065} showed that the map $u^{(\ell)}$ satisfies the identity
\begin{align}
\label{ortho}
   \sum_{j=1}^m y_j\cdot\nabla_ju^{(\ell)}=0.
\end{align}

The following Lemma, which is Lemma\,2.3 from \cite{BS25}, gives an explicit expression for $\Delta ^{k}u^{(\ell)}$ where $k\in\mathbb{N}$.

\begin{Lem}
\label{lem:delta-nak}
 For each $k\in\mathbb{N}$ the map \(u^{(\ell)}\) defined in (\ref{nak-maps}) satisfies 
 \begin{align*}
          \Delta ^{k}u^{(\ell)}=
          \prod_{j=1}^k(2j+\ell-2)(2j-\ell-m)\frac{u^{(\ell)}}{r^{2 k}}.
    \end{align*} 
\end{Lem}

\begin{Bem}
In \cite{MR4817500}, equation (1.6), the following expression for applying the Laplacian \(k\)-times to the radial projection map was presented:
\begin{align*}
\Delta^k\big(\frac{x}{r}\big)=(-1)^k\prod_{j=1}^k(m-2j+1)(2j-1)\frac{x}{r^{2k}}.    
\end{align*}
Clearly, the statement of Lemma \ref{lem:delta-nak} for \(\ell=1\)
is equivalent to the above expression.
\end{Bem}

Now, we establish the following result:

\begin{Satz}
\label{thm:sobolev}
Let $\ell,m\in\mathbb{N}$ with $\ell\leq m$.
The generalized equator map \(u_\star^{(\ell)}\colon B^m\to\s^{{m^\ell}}\) defined
in \eqref{equator-generalized} belongs to \(W^{k,2}(B^m,\s^{{m^\ell}})\) if and only if \(m\geq 2k+1\). In particular, \(u_\star^{(\ell)}\) is an extrinsic $k$-harmonic map.
\end{Satz}

\begin{proof}
Using Lemma\,\ref{lem:delta-nak} it is straightforward to see that
\(u^{(\ell)}\) solves the equation for extrinsic $k$-harmonic maps \eqref{eq:extrinsic-polyharmonic}.
The claim now follows from the fact that the integral 
\begin{align*}
\int_{B^m}\frac{1}{r^q}dx=\vol{(\s^{m-1})}\int_0^1 r^ {m-q}dr,
\end{align*}
where $q\in\mathbb{N}$, is well-defined if and only if $m>q$.
\end{proof}

\section{Stability of extrinsic polyharmonic maps}
\label{sec-stability}
In this section we investigate the stability of the generalized equator map \eqref{equator} considered as a critical point of \eqref{poly-energies}. We prove that understanding the stability of this map reduces to determining the sign of a certain polynomial.

\smallskip

We recall the definitions of stability and energy minimizing, see e.g. \cite{MR4817500}.

\begin{Dfn}
\begin{enumerate}
\item
An extrinsic $k$-harmonic map \(u\colon B^m\to\s^n\) is stable if
\begin{align*}
\frac{d^2}{dt^2}\big|_{t=0}E^{ext}_k(u_t)\geq 0    
\end{align*}
holds for all variations \(u_t\) that satisfy \(u_t-u\in W^{k,2}_0(B^m,\R^{n+1})\). An extrinsic $k$-harmonic map which is not stable, is called unstable.
\item An extrinsic $k$-harmonic map \(u\colon B^m\to\s^n\) is 
called energy minimizing if
\begin{align*}
E^{ext}_k(u)\leq E^{ext}_k(v)  
\end{align*}
holds for all $v\in W^{k,2}(B^m,\s^{n})$ such that \(u-v\in W^{k,2}_0(B^m,\R^{n+1})\).
\end{enumerate}
\end{Dfn}
Clearly, an energy minimizing map is stable.

\smallskip

To characterize the stability of the extrinsic $k$-harmonic map constructed in Theorem\,\ref{thm:sobolev}, we introduce the following polynomials:
\begin{align}
\label{dfn:Q}
Q_1^\ell(m):=&\frac{(m-2)^2}{4}+\ell(2-\ell-m);\\    
Q^\ell_{2s}(m):=&\frac{1}{2^{4s}}\prod_{j=1}^s(m-4j)^2(m+4j-4)^2
-\prod_{j=1}^{2s}(2j+\ell-2)(2j-\ell-m); \\
\nonumber Q^\ell_{2s+1}(m):=&
\frac{(m-2)^2}{2^{4s+2}}\prod_{i=1}^s
(m-4i-2)^2(m+4i-2)^2
+\prod_{j=1}^{2s+1}(2j+\ell-2)(2j-\ell-m),
\end{align}
where $s\in\mathbb{N}$.
Although the polynomials $Q^\ell_{k}(m)$ are defined for $\ell,m\in\mathbb{N}$ only, below we will extend them to $\ell,m
\in\mathbb{R}$ 
whenever necessary. Note that within this interpretation $Q^\ell_{k}(m)$ are continuous functions of $\ell\in\mathbb{R}$ and $m\in\mathbb{R}$.

\smallskip

Our first main result is the following:

\begin{Satz}
\label{thm:stability}
Let \(Q_k^\ell(m)\) be the polynomial defined in \eqref{dfn:Q} and assume that \(m\geq 2k+1\). Then, the generalized equator map \(u_\star^{(\ell)}\colon B^m\to\s^{m^\ell}\) is
\begin{enumerate}
    \item an energy minimizing extrinsic \(k\)-harmonic map if \(Q_k^\ell(m)\geq 0\),
    \item an unstable extrinsic \(k\)-harmonic map if \(Q_k^\ell(m)<0\).
\end{enumerate}
\end{Satz}

Hence, understanding the stability of the generalized equator map 
\(u_\star^{(\ell)}\) boils down to determine for which values of \(m,\ell\)
the polynomial \(Q_k^\ell(m)\) is positive or negative.
Note that by varying the parameter \(\ell\) we can shift the dimension
in which the extrinsic polyharmonic map \(u_\star^{(\ell)}\) turns from being energy minimizing into being unstable. Roughly speaking we find that the larger the value of \(\ell\) the larger is the set of dimensions in which
\(u_\star^{(\ell)}\) is unstable. 

\smallskip

From Theorem\,\ref{thm:stability} we easily obtain the following result for harmonic maps, where we make use of the fact that the sign of $Q_1^\ell(m)=\frac{(m-2)^2}{4}+\ell(2-\ell-m)$ can be determined by solving a quadratic equation.

\begin{Cor}
\label{cor:harmonic}
The generalized equator map \(u_\star^{(\ell)}\colon B^m\to\s^{{m^\ell}}\) defined
in \eqref{equator-generalized} is 
\begin{enumerate}
    \item an energy minimizing harmonic map if
    \begin{align*}
m\geq 2(l+1)+2\sqrt{2}\ell;    
\end{align*}
\item an unstable harmonic map if 
\begin{align*}
    3\leq m <2(l+1)+2\sqrt{2}\ell.
\end{align*}
\end{enumerate}
\end{Cor}


\smallskip

We will postpone the further discussion of the sign of \(Q_k^\ell(m)\) to Section\,\ref{sec-sign} and now turn to the proof of Theorem\,\ref{thm:stability}. We first establish a few preparatory results.

\smallskip

The following Proposition, which was already established in \cite[Proposition 2.2]{MR4436204}, provides an integral inequality for stable critical points 
    of the extrinsic $k$-energy.
\begin{Prop}
\begin{enumerate}
    \item Let \(k=2s\). If \(u\colon B^m\to\s^n\) is a stable critical point 
    of \(E^{ext}_k(u)\) then it satisfies
    \begin{align}
       \int_{B^m}\big(|\Delta^s\eta|^2-\langle\Delta^{2s}u,u\rangle|\eta|^2\big) dx\geq 0 
    \end{align}
     for all \(\eta\in C^\infty_0(B^m,\R^{n+1})\).
     \item Let \(k=2s+1\). If \(u\colon B^m\to\s^n\) is a stable critical point of \(E^{ext}_k(u)\) then it satisfies
    \begin{align}
       \int_{B^m}\big(|\nabla\Delta^s\eta|^2+\langle\Delta^{2s+1}u,u\rangle|\eta|^2\big)dx\geq 0 
    \end{align}
     for all \(\eta\in C^\infty_0(B^m,\R^{n+1})\).
\end{enumerate}    
\end{Prop}

The following Lemma, which was established in \cite[Lemma 2.3]{MR4436204}, provides a sufficient condition for a map to be an energy minimizer of the extrinsic $k$-energy.
\begin{Lem}
\label{lem:energy-minimizing}
Suppose that \(m\geq 2k+1\).
\begin{enumerate}
    \item Let \(k=2s\). If 
    \begin{align}
       \int_{B^m}\big(|\Delta^s\eta|^2-\langle\Delta^{2s}u,u\rangle|\eta|^2\big)dx\geq 0 
    \end{align}
     for all \(\eta\in C^\infty_0(B^m,\R^{n+1})\), then \(u\) is energy minimizing for the extrinsic $k$-energy.
     \item Let \(k=2s+1\). If 
    \begin{align}
       \int_{B^m}\big(|\nabla\Delta^s\eta|^2+\langle\Delta^{2s+1}u,u\rangle|\eta|^2\big)dx\geq 0 
    \end{align}
     for all \(\eta\in C^\infty_0(B^m,\R^{n+1})\), then \(u\) is energy minimizing for the extrinsic $k$-energy.
\end{enumerate}    
\end{Lem}

Further, we recall the following higher order Hardy inequalities, see
\cite{MR2215561} for more background. For this purpose we set, as in \cite{MR4436204},
\begin{align}
\label{constants-hardy}
\alpha_{2s}(m)&=\frac{1}{2^{4s}}\prod_{i=1}^s(m-4i)^2(m+4i-4)^2, \\
\nonumber\alpha_{2s+1}(m)&=\frac{1}{2^{4s+2}}(m-2)^2\prod_{i=1}^s
(m-4i-2)^2(m+4i-2)^2,\\
\alpha_1(m)\nonumber&=\frac{(m-2)^2}{4},
\end{align}
where \(s\geq 1\).

\begin{Satz}[Hardy inequalities]
\label{thm:hardy}
Let \(\alpha_k(m)\) be the constants defined in \eqref{constants-hardy}.
\begin{enumerate}
    \item Let \(k=2s\) and \(m\geq 2k+1\). Then, we have
     \begin{align*}
        \int_{B^m}|\Delta^s\eta|^2dx\geq\alpha_{2s}(m)
        \int_{B^m}\frac{|\eta|^2}{r^{4s}}dx,\qquad s\geq 1 
     \end{align*}
     for all \(\eta\in W^{k,2}_0(B^m,\R^{n+1})\).
    \item Let \(k=2s+1\) and \(m\geq 2k+1\). Then, we have
     \begin{align*}
        \int_{B^m}|\nabla\Delta^s\eta|^2dx\geq\alpha_{2s+1}(m)
        \int_{B^m}\frac{|\eta|^2}{r^{4s+2}}dx,\qquad s\geq 0
     \end{align*}
     for all \(\eta\in W^{k,2}_0(B^m,\R^{n+1})\).
\end{enumerate}
\end{Satz}

\begin{proof}
The above Hardy inequalities can be found in \cite{MR2215561} for functions on \(B^m\). However, using the arguments presented in 
\cite[p. 272]{MR4817500} these inequalities also hold in the case for functions taking values in \(\R^{n+1}\).
\end{proof}

Combining the previous results we now prove Theorem\,\ref{thm:stability}:
\begin{proof}[Proof of Theorem \ref{thm:stability}]
The proof is analogous to that of  \cite[Theorem 1.1]{MR4436204}.

\smallskip

We prove (1) for the even case, i.e. for $k=2s$. The case $k=2s+1$ is proven analogously.
We set
\begin{align*}
    B_k^{\ell}(m):=\prod_{j=1}^k(2j+\ell-2)(2j-\ell-m).
\end{align*}
Using Lemma\,\ref{lem:delta-nak} we thus get
\begin{align}
\label{eq:delta2s}
\langle \Delta^{2s}u^{(\ell)},u^{(\ell)}\rangle=B_{2s}^{\ell}(m)r^{-4s},  
\end{align}
which also holds for \(u_\star^{(\ell)}\)
Moreover, using (\ref{eq:delta2s}) and Theorem\,\ref{thm:hardy} (1), we find
\begin{align*}
 \int_{B^m}\big(\lvert\Delta^s\phi\rvert^2-\langle \Delta^{2s}u_\star^{(\ell)},u_\star^{(\ell)}\rangle\lvert\phi\rvert^2\big)dx\geq Q_{2s}^{\ell}(m)\int_{B^m}\lvert\phi\rvert^2r^{-4s}dx.   
\end{align*}
Since $Q_{2s}^{\ell}(m)\geq 0$ by assumption, 
we get
\begin{align*}
 \int_{B^m}\big(\lvert\Delta^s\phi\rvert^2-\langle \Delta^{2s}u_\star^{(\ell)},u_\star^{(\ell)}\rangle\lvert\phi\rvert^2\big)dx\geq 0,  
\end{align*}
hence the claim follows from Lemma\,\ref{lem:energy-minimizing}, (1). 

\smallskip

We now prove (2) for $k=2s$. The case $k=2s+1$ is proven analogously.
Let $\epsilon>0$ be such that $Q_{2s}^{\ell}(m)+\epsilon<0$.
The density of $C^{\infty}_0(B^m,\mathbb{R}^{m^{\ell}+1})$ in $W_{0}^{k,2}(B^m,\mathbb{R}^{m^{\ell}+1})$ yields
\begin{align*}
\alpha_{2s}(m)=\mbox{inf}\{\frac{\int_{B^m}\lvert\Delta^{s}\phi\rvert^2dx}{\int_{B^m}\frac{\lvert\phi\rvert^2}{r^{4s}}dx}\,\lvert\,\phi\in C^{\infty}_0(B^m,\mathbb{R}^{m^{\ell}+1})\},    
\end{align*}
see \cite[page 9164]{MR4436204}. Consequently, there exists a $\phi_{\epsilon}\in C^{\infty}_0(B^m,\mathbb{R}^{m^{\ell}+1})$ such that
\begin{align*}
\frac{\int_{B^m}\lvert\Delta^{s}\phi_{\epsilon}\rvert^2dx}{\int_{B^m}\frac{\lvert\phi_{\epsilon}\rvert^2}{r^{4s}}dx}\leq\alpha_{2s}(m)+\epsilon.
\end{align*}
Using this inequality and (\ref{eq:delta2s}) we obtain
\begin{align*}
 \int_{B^m}\big(\lvert\Delta^s\phi_{\epsilon}\rvert^2-\langle \Delta^{2s}u_\star^{(\ell)},u_\star^{(\ell)}\rangle\lvert\phi_{\epsilon}\rvert^2\big)dx\leq (Q_{2s}^{\ell}(m)+\epsilon)\int_{B^m}\lvert\phi_{\epsilon}\rvert^2r^{-4s}dx<0,      
\end{align*}
where we used the assumption of the Theorem to deduce the last inequality. Hence the claim is established.
\end{proof}

\section{Analysis of the sign of \texorpdfstring{$Q_k^{\ell}(m)$}{Q}}
\label{sec-sign}
As discussed in Section\,\ref{sec-stability} after the statement of Theorem\,\ref{thm:stability},
understanding the stability of the generalized equator map 
\(u_\star^{(\ell)}\) boils down to determine for which values of \(m,\ell\)
the polynomial \(Q_k^\ell(m)\) is positive or negative.
We will be concerned with this problem in the current section
and present a detailed analysis.
Note that the first result of this kind is already contained in the previous section, namely Corollary\,\ref{cor:harmonic}.

\smallskip

We first show that  $Q^\ell_{k}(2k+1)<0$
for all \(k\geq 1\) and \(1\leq\ell\leq 2k+1\).
Combining this statement with Theorem\,\ref{thm:stability} yields that for specific choices of $\ell$ the generalized equator map is unstable, see Corollary\,\ref{cor-unstable} for the precise details.

\begin{Lem}
\label{negative}
We have
\begin{align}
  Q^\ell_{k}(2k+1)<0
\end{align}
for all \(k\geq 1\) and \(1\leq\ell\leq 2k+1\).
\end{Lem}

\begin{proof}
We provide the proof for \(k=2s\). The case \(s=2k+1\) is proven analogously.
To this end we set
\begin{align*}
Q^\ell_{2s}(m):=&\underbrace{\frac{1}{16^s}\prod_{j=1}^s(m-4j)^2(m+4j-4)^2}_{:=A_{2s}(m)}
-\underbrace{\prod_{j=1}^{2s}(2j+\ell-2)(2j-\ell-m)}_{:=B^\ell_{2s}(m)}.    
\end{align*}
Concerning the first term we have the identity
\begin{align*}
A_{2s}(2k+1)= \prod_{j=0}^{2s-1}\big(4j^2+2j+\frac{1}{4}\big),
\end{align*}
see \cite[Proof of Lemma 1.7]{MR4436204}. Concerning the second term we obtain
\begin{align*}
B^\ell_{2s}(4s+1)=&\prod_{j=1}^{2s}(2j+\ell-2)(2j-\ell-4s-1) =\prod_{j=0}^{2s-1}(2j+\ell)\prod_{j=0}^{2s-1}(2j-\ell-4s+1) \\
=&\prod_{j=0}^{2s-1}(2j+\ell)\prod_{j=0}^{2s-1}(1+2j+\ell)=\prod_{j=0}^{2s-1}(4j^2+\ell^2+2j+4jl+\ell).
\end{align*}
Using these two identities we find 
\begin{align*}
Q^\ell_{2s}(4s+1)=\prod_{j=0}^{2s-1}\big(4j^2+2j+\frac{1}{4}\big)
-\prod_{j=0}^{2s-1}(4j^2+\ell^2+2j+4jl+\ell)
<0,
\end{align*}
completing the proof.
\end{proof}

\begin{Cor}
\label{cor-unstable}
   The generalized equator map $u_\star^{(\ell)}:B^{2k+1}\rightarrow\s^{{(2k+1)}^{\ell}}$ is an unstable extrinsic $k$-harmonic map for all \(k\geq 1\) and \(1\leq\ell\leq 2k+1\).
\end{Cor}

\smallskip

Next we prove that for given \(k,\ell\geq 1\) there exists an integer \(m^\ast=m(k,\ell)>2k+1\) such that $Q_k^\ell(m)$ is non-negative if \(m\geq m^\ast\) and negative if \(2k+1\leq m<m^\ast\).

\begin{Lem}
\label{lem:m}
For all \(k\geq 1\) there exists an integer \(m^\ast=m(k,\ell)>2k+1\) such that
\begin{enumerate}
    \item \(Q_k^\ell(m)\geq 0\) if \(m\geq m^\ast\); 
    \item \(Q_k(m)<0\) if \(2k+1\leq m<m^\ast\).
\end{enumerate}
\end{Lem}
\begin{proof}
The proof is inspired by that of \cite[Lemma 1.7]{MR4436204}.
We will give the proof for $k=2s$, the case $k=2s+1$ is proven analogously.

\smallskip

We start by rewriting $Q^\ell_{2s}(m)$:
\begin{align}
\label{q-r}
Q^\ell_{2s}(m)=\prod_{j=1}^{2s}(\frac{m}{2}+2(j-1)-2s)^2
-\prod_{j=1}^{2s}(\ell+2j-2)(m+\ell-4s+2j-2).
\end{align}
This identity indeed holds: From the very definition of $Q^\ell_{2s}(m)$
we get
\begin{align*}
Q^\ell_{2s}(m)=\frac{1}{2^{4s}}\prod_{j=1}^s(m-4j)^2(m+4j-4)^2
-\prod_{j=1}^{2s}(\ell+2j-2)(m+\ell-2j).
\end{align*}
We have 
\begin{align*}
    \prod_{j=1}^s(m-4j)^2(m+4j-4)^2
    = \prod_{k=-(s-1)}^s(m-4k)^2
    =\prod_{j=1}^{2s}(m+4(j-1)-4s)^2,
\end{align*}
where for the last equality we shifted the indices, more precisely we used the substitution $k=s-(j-1)$. 
Further, since 
\begin{align*}
 \prod_{j=1}^{2s}(m+\ell-2j)=\prod_{j=1}^{2s}(m+\ell-4s+2j-2),   
\end{align*}
identity (\ref{q-r}) holds.

\smallskip

We now show the existence of $m^\ast=m(k,\ell)$.
By Lemma\,\ref{negative} we have 
\begin{align}
  Q^\ell_{k}(2k+1)<0
\end{align}
for all \(k\geq 1\) and \(1\leq\ell\leq 2k+1\). 
Further, the highest order coefficient in $m$ of $Q^\ell_{2s}(m)$ is positive, hence $\lim_{m\rightarrow\infty}Q^\ell_{2s}(m)=\infty$. Consequently, there exists an integer $N>4s+1$ such that $Q^\ell_{2s}(N)>0$.

\smallskip

We complete the proof by an induction argument, i.e. by showing that if there exists $N>4s+1$, $N\in\mathbb{N}$, with $Q_{2s}^{\ell}(N)\geq 0$, then $Q_{2s}^{\ell}(N+1)>0$.

Note that the equation
\begin{align}
\label{q-gb}
 Q_{2s}^{\ell}(N+1)=\prod_{j=1}^{2s}\gamma_j\prod_{j=1}^{2s}(\frac{N}{2}+2(j-1)-2s)^2-\prod_{j=1}^{2s}\beta_j\prod_{j=1}^{2s}(\ell+2j-2)(N+\ell-4s+2j-2)   
\end{align}
with 
\begin{align*}
\gamma_j=\bigg(\frac{\frac{N+1}{2}+2(j-1)-2s}{\frac{N}{2}+2(j-1)-2s}\bigg)^2    
\end{align*}
and
\begin{align*}
 \beta_j=\frac{N+1+\ell-4s+2j-2}{N+\ell-4s+2j-2},   
\end{align*}
$j\in{1,\dots,2s}$,
holds.
We rewrite (\ref{q-gb}) as
\begin{align*}
 Q_{2s}^{\ell}(N+1)=(\prod_{j=1}^{2s}\gamma_j)\,Q_{2s}^{\ell}(N)+(\prod_{j=1}^{2s}\gamma_j-\prod_{j=1}^{2s}\beta_j)\,\prod_{j=1}^{2s}(\ell+2j-2)(N+\ell-4s+2j-2).   
\end{align*}
Write $N=4s+1+a$, with $a>0$.
Then
\begin{align*}
\gamma_j-\beta_j=\frac{-4+a^2+6j+a(-1+2\ell+4j)+ 
 \ell(-5+8j)}{(-1 + a + \ell+ 2 j) (-3 + a + 4 j)^2},
\end{align*}
which is positive since $j\geq 1$.
Thus 
\begin{align*}
\prod_{j=1}^{2s}\gamma_j-\prod_{j=1}^{2s}\beta_j>0   
\end{align*}
and therefore $ Q_{2s}^{\ell}(N+1)>0$, whence the claim.
\end{proof}

For $k=2$ and $\ell=1$ we can easily determine that $m(2,1)=10$ and hence
Lemma\,\ref{lem:m} recovers the result of Hong and Thompson \cite{MR2322746} for $k=2$ and $\ell=1$.

\smallskip

Next, we attack the problem of determining information about $m^\ast=m(k,\ell)$.
With the help of a computer algebra system we obtain the following values for $m(k,\ell)$
with $k\leq 4$ and $\ell\leq 10$.

\begin{center}
\captionof{table}{$m(k,\ell)$ for $k\leq 4$ and $\ell\leq 10$}
\begin{tabular}{ |c||c|c|c|c|c|c|c|c|c|c| } 
 \hline
$\ell$ & 1 & 2 & 3 & 4& 5 & 6& 7 & 8& 9 & 10\\ 
 \hline
  \hline
$m(1,\ell)$ & 7 & 12& 17 & 22& 27 & 31& 36 & 41& 46 & 51 \\ 
 \hline
$m(2,\ell)$ & 10 & 15& 20 & 25& 30 & 35& 39 & 44& 49 & 54 \\ 
 \hline
$m(3,\ell)$ & 12 & 18& 23 & 28& 33 & 38& 43 & 47& 52 & 57 \\  \hline
$m(4,\ell)$ & 15 & 20& 25 & 31& 36 & 41& 46 & 50& 55 & 60 \\ 
 \hline
\end{tabular}
\end{center}

\medskip

Note that $Q^{\ell+1}_{k}(m)-Q^{\ell}_{k}(m)< 0$ for any $k\in\mathbb{N}$. Hence $m(k,\ell+1)> m(k,\ell)$ for any $k,\ell\in\mathbb{N}$. For fixed $k$, when $\ell$ increases, the dimension $m\geq \ell$ for which $Q^{\ell}_{k}(m)< 0$ increases as well. Roughly speaking, we lose stability when increasing $\ell$ and we can make a transition to instability as quickly as desired.

\smallskip

For $k=2$ and $k=3$ the equation $Q^{\ell}_{k}(m)=0$ can be solved explicitly for $m$. However, since the expressions are rather long and do not contribute to further insight, we do not state them here. Instead, we will now focus on the case $k=2$ in which it is possible to solve the equation $Q^{\ell}_{k}(m)=0$ explicitly for $\ell$ which gives a good impression of the general picture.

\begin{Lem}
 Let $m\in\mathbb{N}$ be given.
 Set $$L:=1-\frac{m}{2}+\frac{1}{2}\sqrt{20-8m+m^2+(m-4)\sqrt{m^2+16}}.$$
 Then, we have
 \begin{enumerate}
     \item $Q^\ell_2(m)\geq 0$ for any $\ell\in[1,L]$;
     \item $Q^\ell_2(m)< 0$ for any $\ell\in (L,m]$.
 \end{enumerate}
\end{Lem}
\begin{proof}
Fix $m\in\mathbb{N}$.
From (\ref{dfn:Q}) we get
\begin{align*}
Q^\ell_2(m)=\frac{1}{16}m^2(m-4)^2-\ell(\ell+2)(m+\ell-2)(m+\ell-4).
\end{align*}
    A straightforward computation gives that $Q^\ell_2(m)=0$ if and only if 
\begin{align}
\label{eq:l}
 \ell=1-\frac{m}{2}\pm\frac{1}{2}\sqrt{20-8m+m^2\pm(m-4)\sqrt{m^2+16}}. 
\end{align}
Furthermore, we have 
\begin{align*}
L&\leq 1-\frac{m}{2}+\frac{1}{2}\sqrt{20-8m+m^2+(m-4)(m+4)}= 1-\frac{m}{2}+\frac{1}{2}\sqrt{2}\sqrt{m^2-4m+2}\\
&\leq  1-\frac{m}{2}+\frac{1}{2}\sqrt{2}(m-2)=(\sqrt{2}-1)(\frac{m}{2}-1)<m.
\end{align*}
Further $Q^m_2(m)=-16m+17m^2+\frac{7}{2}m^3-\frac{63}{16}m^4<0$ for any $m\geq 5$.
Hence, and since $\ell\leq L$ for any $\ell$ satisfying (\ref{eq:l}),
continuity of $Q^\ell_2(m)$ with respect to $\ell$ yields
$Q^\ell_2(m)<0$ for any $\ell\in(L,m]$.

\smallskip

We have $Q^\ell_2(m)\geq 0$ for any $m\in\mathbb{N}$ given and any $\ell\in[1,L]$.
Indeed, 
observe that the only other possible non-negative $\ell\in\mathbb{R}$ in \eqref{eq:l} is given by
\begin{align*}
 L_1=1-\frac{m}{2}+\frac{1}{2}\sqrt{20-8m+m^2-(m-4)\sqrt{m^2+16}}.  
\end{align*}
However,
\begin{align*}
 L_1&=1-\frac{m}{2}+\frac{1}{2}\sqrt{20-8m+m^2-(m-4)\sqrt{m^2+16}}\leq 1-\frac{m}{2}+\frac{1}{2}(m-4)=-1,
\end{align*}
since $m\geq 5$ and thus $(m-4)\sqrt{m^2+16}\geq 4$. Consequently, the polynomial $Q^\ell_2(m)$ does not change its sign when $\ell$ varies in the interval $[1,L)$.
\end{proof}

Note that from the very definition of (\ref{dfn:Q}) we obtain that 
for $\epsilon>0$ we have $Q^\ell_k(m)>0$ for all $m\geq (2(1+\sqrt{2})+\epsilon)\ell$ whenever $\ell$ is sufficiently large.
Clearly, this gives that $m(k,\ell)$ is bounded from above, however this bound is not explicit.
Below we provide an explicit upper bound for $m(k,\ell)$.
In the next Lemma we do so for the specific case $k=2$.

\begin{Lem}
\label{k2}
We have $Q^\ell_2(m)>0$ for all $m\geq 5(\ell+1)$. 
Therefore $m(2,\ell)\leq 5(\ell+1)$ for any $\ell\in\mathbb{N}$.   
\end{Lem}
\begin{proof}
  Evaluating $Q^\ell_2(m)$
  for $m=5(\ell+1)$
  yields
  $$Q^\ell_2(5(\ell+1))=\frac{1}{16}(25+204\ell+334\ell^2-36\ell^3+49\ell^4),$$ which is positive for any $\ell\in\mathbb{N}$.
  The claim thus follows from Lemma\,\ref{lem:m}.
\end{proof}  
 
Our next goal is to provide an upper bound for $m(k,\ell)$.
For this purpose we give a preparatory lemma, whose proof is contained in Appendix\,A:

\begin{Lem}
\label{app}
Let $\ell,s\in\mathbb{N}$ and introduce $f_1,f_2:\mathbb{N}\times\mathbb{N}\mapsto\mathbb{R}$ by
\begin{align*}
f_1(\ell,s)&=(3\ell-\frac{1}{2}+s)(3\ell+\frac{1}{2}+s)\prod_{i=0}^4(\frac{5}{4}\ell-\frac{1}{4}+\frac{5}{2}s+i)^2,\\
f_2(\ell,s)&=(\frac{5}{4}\ell-\frac{1}{4}+\frac{s}{2})^2\prod_{j=0}^3(\frac{\ell}{2}+2s+j)\prod_{i=0}^5(3\ell-\frac{1}{2}+3s+i).
\end{align*}
We have that $f_1-f_2>0$.
\end{Lem}

With the help of 
the following Theorem we get an upper bound for $m(k,\ell)$, see Corollary\,\ref{cor-main}.

\begin{Thm}
\label{con}
For any $k\geq 2$ and $\ell\in\mathbb{N}$ we have $Q_k^{\ell}(m)>0$ for any $m\geq 5+3(k-2)+5\ell$. If this inequality is satisfied, the generalized equator map $u_\star^{(\ell)}:B^m\rightarrow\s^{m^\ell}$ is an energy minimizing extrinsic $k$-harmonic map.
\end{Thm}

\begin{proof}
We prove the result for $k$ even, the proof for $k$ odd works analogously, making use of induction on $k$.

\smallskip

For $x\in\mathbb{R}_+$ and $n\in\mathbb{N}$ we have
\begin{align*}
    x(x+1)\dots(x+n-1)=\frac{\Gamma(x+n)}{\Gamma(x)}, 
\end{align*}
where $\Gamma:\mathbb{R}_+\rightarrow\mathbb{R}$ denotes the Gamma function.
Thus we obtain
\begin{align*}
    x(x+2)\dots(x+2(n-1))=2^n\frac{\Gamma(\frac{x}{2}+n)}{\Gamma(\frac{x}{2})}. 
\end{align*}
Consequently, we get
    \begin{align}
    \label{q-gamma}
          Q_{2s}^{\ell}(5+3(2s-2)+5\ell)=2^{4s}(\frac{\Gamma(\frac{5}{4}\ell-\frac{1}{4}+\frac{5}{2}s)^2}{\Gamma(\frac{5}{4}\ell-\frac{1}{4}+\frac{1}{2}s)^2}-\frac{\Gamma(3\ell-\frac{1}{2}+3s)\Gamma(\frac{1}{2}\ell+2s)}{\Gamma(3\ell-\frac{1}{2}+s)\Gamma(\frac{1}{2}\ell)}).
    \end{align}

We will first show $Q_{2s}^{\ell}(5+3(2s-2)+5\ell)>0$ for $\ell,s\geq 1$. Since $\ell,s\geq 1$, all arguments of the Gamma function appearing in (\ref{q-gamma}) are  positive and thus we have that  $Q_{2s}^{\ell}(5+3(2s-2)+5\ell)>0$ if and only if 
\begin{align}
\label{eqeq}
\Gamma(\frac{5}{4}\ell-\frac{1}{4}+\frac{5}{2}s)^2\Gamma(3\ell-\frac{1}{2}+s)\Gamma(\frac{1}{2}\ell)-\Gamma(\frac{5}{4}\ell-\frac{1}{4}+\frac{1}{2}s)^2\Gamma(3\ell-\frac{1}{2}+3s)\Gamma(\frac{1}{2}\ell+2s)>0
\end{align}
is satisfied.

\smallskip

We will now prove by induction that inequality (\ref{q-gamma}) holds true for all $s\in\mathbb{N}$.
More precisely, we show that if this inequality holds true for $s_0\in\mathbb{N}$, then it also holds true for $s_0+2$.
Further, we prove inequality (\ref{q-gamma}) for $s_0=1$ and $s_0=2$. Overall, this gives us inequality (\ref{q-gamma}) for all $s\in\mathbb{N}$.

The reason for this seemingly cumbersome procedure is, when substituting $s$ by $s+2$, all arguments of the Gamma functions appearing in (\ref{eqeq}) are shifted only by integers
such that we can apply
the functional equation $\Gamma(x+1)=x\Gamma(x)$, $x\in\mathbb{R}$, of the Gamma function in a later step.

\smallskip

The induction basis is as follows:
The claim has been verified for $k=2$ (i.e. $s=1$) in Lemma\,\ref{k2}.
For $k=4$ (i.e. $s=2$) we have 
\begin{align*}
  Q_{4}^{\ell}(11+5\ell)=\frac{1}{256}(&12006225+64479240\ell+ 132331788\ell^2 + 131725656\ell^3 + 
   66283606\ell^4\\ &+ 15832440\ell^5 + 1765612\ell^6 + 316584\ell^7 + 58849\ell^8),   
\end{align*}
which is positive for any $\ell\in\mathbb{N}$.
  The claim for $k=4$ thus follows from Lemma\,\ref{lem:m}.

\smallskip

We now provide the induction step, which is the same for $s$ even and $s$ odd, so here we do not make a distinction between these cases.
Assume $Q_{2s}^{\ell}(5+3(2s-2)+5\ell)>0$ for a fixed $s\in\mathbb{N}$ and any $\ell\in\mathbb{N}$. We will prove $Q_{2(s+2)}^{\ell}(5+3(2(s+2)-2)+5\ell)>0$. From (\ref{eqeq}) we get that this inequality is equivalent to 
\begin{multline*}
f_1(\ell,s)\Gamma(\frac{5}{4}\ell-\frac{1}{4}+\frac{5}{2}s)^2\Gamma(3\ell-\frac{1}{2}+s)\Gamma(\frac{1}{2}\ell)\\ -f_2(\ell,s)\Gamma(\frac{5}{4}\ell-\frac{1}{4}+\frac{1}{2}s)^2\Gamma(3\ell-\frac{1}{2}+3s)\Gamma(\frac{1}{2}\ell+2s)>0,
\end{multline*}
where $f_1$ and $f_2$ are as in Lemma\,\ref{app}.
Note that $f_1>0$ and $f_2>0$.
We have $f_1-f_2>0$, see Appendix\,A.
Thus we get
from the induction assumption, i.e. that (\ref{eqeq}) is satisfied, the inequality
\begin{align*}
&f_1(\ell,s)\Gamma(\frac{5}{4}\ell-\frac{1}{4}+\frac{5}{2}s)^2\Gamma(3\ell-\frac{1}{2}+s)\Gamma(\frac{1}{2}\ell) \\
&-f_2(\ell,s)\Gamma(\frac{5}{4}\ell-\frac{1}{4}+\frac{1}{2}s)^2\Gamma(3\ell-\frac{1}{2}+3s)\Gamma(\frac{1}{2}\ell+2s)\\
>&f_1(\ell,s)\big(\Gamma(\frac{5}{4}\ell-\frac{1}{4}+\frac{5}{2}s)^2\Gamma(3\ell-\frac{1}{2}+s)\Gamma(\frac{1}{2}\ell)-\Gamma(\frac{5}{4}\ell-\frac{1}{4}+\frac{1}{2}s)^2\Gamma(3\ell-\frac{1}{2}+3s)\Gamma(\frac{1}{2}\ell+2s)\big) \\
>&0.
\end{align*}
Hence $Q_{2(s+2)}^{\ell}(5+3(2(s+2)-2)+5\ell)>0$.
 The claim of the theorem now follows from Lemma\,\ref{lem:m}.
\end{proof}

As an immediate consequence we obtain the following upper bound for $m(k,\ell)$:

\begin{Cor}
\label{cor-main}
We have $m(k,\ell)\leq 5+3(k-2)+5\ell$ for any $k\geq 2$.
\end{Cor}

\smallskip

In the following Theorem we show that if \(\ell=m\geq 2k+1\), the generalized equator map \(u_\star^{(m)}\colon B^m\rightarrow\s^{m^m}\) is an unstable extrinsic $k$-harmonic map.    

\begin{Satz}
\label{thm:lm}
Let \(Q_k^\ell(m)\) be the polynomial defined in \eqref{dfn:Q}. Assume that \(\ell=m\geq 2k+1\). Then, the generalized equator map \(u_\star^{(m)}\colon B^m\to\s^{m^m}\) is unstable.   
\end{Satz}

\begin{proof}
 We consider the case $k=2s$, $s\in\mathbb{N}$, first.
 From the very definition of $Q_{2s}^{\ell}(m)$ and the proof of Lemma\,\ref{lem:m} we get
\begin{align*}
Q^m_{2s}(m)=\prod_{j=1}^{2s}\big(\frac{m}{2}+2(j-1)-2s\big)^2
-2^{2s}\prod_{j=1}^{2s}(m+2j-2)(m-j).
\end{align*}
Observe that for any $s\in\mathbb{N}$ and any $m\geq 4s+1$ we have 
\begin{align*}
\frac{m}{2}+2(j-1)-2s\geq 0,\quad m+2j-2\geq 0,\quad m-j\geq 0,
\end{align*}
for each $j\in\{1,\dots,2s\}.$
Furthermore, we have
\begin{align*}
\frac{m}{2}+2(j-1)-2s\leq m+2j-2   
\end{align*}
and also the inequality
\begin{align*}
\frac{m}{2}+2(j-1)-2s\leq 4(m-j)   
\end{align*}
holds due to $m\geq 4s+1$ and $j\in\{1,\dots,2s\}.$
Consequently, we get
\begin{align*}
Q^m_{2s}(m)=&\prod_{j=1}^{2s}(\frac{m}{2}+2(j-1)-2s)^2
-2^{2s}\prod_{j=1}^{2s}(m+2j-2)(m-j)\\
\leq &\prod_{j=1}^{2s}(\frac{m}{2}+2(j-1)-2s)^2
-2^{2s}\prod_{j=1}^{2s}(\frac{m}{2}+2(j-1)-2s)(m-j)\\
\leq &\prod_{j=1}^{2s}(\frac{m}{2}+2(j-1)-2s)^2
-2^{2s-2}\prod_{j=1}^{2s}(\frac{m}{2}+2(j-1)-2s)^2 \\
<&0,
\end{align*}
whence the claim for the case $k=2s$ follows from Theorem\,\ref{thm:stability}.

\smallskip

Next we consider the case $k=2s+1$, $s\in\mathbb{N}$.
From (\ref{dfn:Q}) we obtain
\begin{align*}
    Q^m_{2s+1}(m)=
\frac{(m-2)^2}{2^{4s+2}}\prod_{i=1}^s
(m-4i-2)^2(m+4i-2)^2
-2^{2s+1}\prod_{j=1}^{2s+1}(m+2j-2)(m-j).
\end{align*}
Since 
\begin{align*}
(m-2)^2\prod_{i=1}^s
(m-4i-2)^2(m+4i-2)^2=\prod_{i=1}^{2s+1}
(m-2-4(i-s-1))^2,    
\end{align*}
we deduce that
\begin{align*}
    Q^m_{2s+1}(m)&=
\frac{1}{2^{4s+2}}\prod_{i=1}^{2s+1}
\big(m+2-4(i-s)\big)^2-2^{2s+1}\prod_{j=1}^{2s+1}(m+2j-2)(m-j)\\
&=\prod_{i=1}^{2s+1}
(\frac{m}{2}+2(s-i)+1)^2-2^{2s+1}\prod_{j=1}^{2s+1}(m+2j-2)(m-j).
\end{align*}
Further, since
\begin{align*}
\frac{m}{2}+2(s-i)+1\leq m+2i-2   
\end{align*}
is equivalent to
\begin{align*}
4s+6-8i\leq m
\end{align*}
and $i\in\{1,\dots,2s+1\}$, $m\geq 4s+3$, the second and thus the first inequality are satisfied.
In addition, the inequality
\begin{align*}
\frac{m}{2}+2(s-i)+1\leq 4(m-i)  
\end{align*}
can be rewritten as
\begin{align*}
\frac{4}{7}s+\frac{4}{7}i+\frac{2}{7}\leq m
\end{align*}
and as $i\in\{1,\dots,2s+1\}$, $m\geq 4s+3$, again the second and thus the first inequality hold.
Hence, we have
\begin{align*}
    Q^m_{2s+1}(m)<0,
\end{align*}
and the claim for $k=2s+1$, $s\in\mathbb{N}$ now follows from Theorem\,\ref{thm:stability}.
\end{proof}

\section{Stability of the generalized equator map as \texorpdfstring{$p$}{p}-harmonic maps}
\label{sec-stability-p}
In this section we will prove the results concerning
the stability of the generalized equator map \eqref{equator-generalized}, 
when considered as a \(p\)-harmonic map, that have already been presented in the introduction.

The map \(u^{(\ell)}\) defined in \eqref{nak-maps} is also $p$-harmonic, that is a critical point of
\eqref{p-energy}, i.e. it solves \eqref{p-tension-field-sphere},
see the appendix of \cite{MR4757446} for the precise details.
Hence, also the generalized equator map \(u_\star^{(\ell)}\) given by \eqref{equator-generalized} is \(p\)-harmonic.

\begin{Satz}
\label{thm:p-harmonic}
Assume that \(m>p\geq 2\). 
\begin{enumerate}
    \item If
     \begin{align*}
     \frac{(m-p)^2}{4\ell(\ell+m-2)} < 1
     \end{align*}
     then the generalized equator map \(u_\star^{(\ell)}\colon B^m\rightarrow\s^{m^\ell}\) is an unstable \(p\)-harmonic map.
     \item If
     \begin{align*}
        \frac{(m-p)^2}{4\ell(\ell+m-2)}\geq  1
     \end{align*}
     then the generalized equator map \(u_\star^{(\ell)}\colon B^m\rightarrow\s^{m^\ell}\) is the unique energy minimizer 
     of the \(p\)-energy \eqref{p-energy}.
\end{enumerate}
\end{Satz}

The proof of this result is inspired from the classic 
results of Jäger and Kaul \cite{MR705882} and the corresponding generalization to the \(p\)-harmonic case of Hong \cite{MR1743430}.

\smallskip
In order to establish Theorem \ref{thm:p-harmonic}
let us recall the second variation formula of the \(p\)-energy
\eqref{p-energy} for maps to the sphere
\begin{align}
\label{eq:sv-p-energy}
\frac{d^2}{dt^2}E_p(u_t)\big|_{t=0}=\int_M|du|^{p-2}\big(|d\eta|^2-|du|^2|\eta|^2\big)\dv+(p-2)\int_M|du|^{p-4}|\langle du,d\eta\rangle|^2\dv.
\end{align}
Here, \(\eta=\frac{du_t}{dt}\big|_{t=0}\) and 
\(\langle\eta,u\rangle=0\) as \(\eta\in\Gamma(T\s^n)\).
For a derivation of 
\eqref{eq:sv-p-energy} we refer to \cite[p. 143]{MR1619840},
see also \cite[Def. 2.1]{MR4715070}.

First, we establish the instability statement of Theorem \ref{thm:p-harmonic}:

\begin{Prop}
If the inequality
\begin{align}
\label{eq:condition-p-harmonic-instability}
\frac{(m-p)^2}{4\ell(\ell+m-2)} < 1
\end{align}
is satisfied, then \(u_\star^{(\ell)}\colon B^m\rightarrow\s^{m^\ell}\) is an unstable \(p\)-harmonic map.  
\end{Prop}

\begin{proof}
We employ the second variation formula \eqref{eq:sv-p-energy}.
We choose \(\eta=v(r)e_{\R^{m^\ell+1}}\) as in 
\cite[Theorem 2]{MR705882}, where $v\in C^2(\mathbb{R})$
 will be determined along the proof.
Since \(\langle\eta,u_\star^{(\ell)}\rangle=0\) and as we have
\(\Delta u_\star^{(\ell)}=-\frac{\ell(\ell+m-2)}{r^2}u_\star^{(\ell)}\),
we get \(\langle du_\star^{(\ell)},d\eta\rangle=-\langle u_\star^{(\ell)},\Delta\eta\rangle\).

Inserting these data into the formula for the second variation \eqref{eq:sv-p-energy}, we obtain
\begin{align}
\label{eq:sv-p-reduced}
\hess E_p(u^{(\ell)})(\eta,\eta)
=&\ell^\frac{p-2}{2}(\ell+m-2)^\frac{p-2}{2}\int_{B^m}
\frac{1}{r^{p-2}}\big(v'(r)^2-\frac{\ell(\ell+m-2)}{r^2}v^2\big)dx\\
\nonumber=&\vol(\s^{m-1})\ell^\frac{p-2}{2}(\ell+m-2)^\frac{p-2}{2}\int_0^1
r^{m-p+1}\big(v'(r)^2-\frac{\ell(\ell+m-2)}{r^2}v^2\big)dr.
\end{align}

We suppose that \(v(r)\) satisfies the following ordinary differential equation
\begin{align}
\label{eq:ode-v-p-harmonic}
\begin{cases}
v''(r)+\frac{m-p+1}{r}v'(r)+\frac{\ell(\ell+m-2)-\epsilon}{r^2}v(r)=0  &\text{ on } r_0\leq r\leq 1, \\
v(r)=0 &\text{ on } r\leq r_0,
\end{cases}    
\end{align}
where \(\epsilon>0\) is a small number.
If we insert the solution of \eqref{eq:ode-v-p-harmonic}
into \eqref{eq:sv-p-reduced}
and after conducting integration by parts we find
\begin{align*}
\hess E_p(u_\star^{(\ell)})(\eta,\eta)
=-\epsilon\vol(\s^{m-1})\ell^\frac{p-2}{2}(\ell+m-2)^\frac{p-2}{2}\int_{r_0}^1v^2(r)r^{m-p-1}dr<0    
\end{align*}
producing a negative sign for the Hessian of the \(p\)-energy.

Hence, it remains to analyze under which conditions we can 
solve the ordinary differential equation \eqref{eq:ode-v-p-harmonic}.
To this end, we perform the transformation \(z(e^t)=v(t)\) such that we get
\begin{align*}
z''+(m-p)z'+\big(\ell(\ell+m-2)-\epsilon\big)z=0.  
\end{align*}
Solving this equation and transforming back we obtain
\begin{align*}
v(r)=
\begin{cases}
    &r^{\frac{p-m}{2}}\big(c_1\cos(\sqrt{-\mu}\ln r)+c_2\sin(\sqrt{-\mu}\ln r)\big), \\
    &r^{\frac{p-m}{2}}\big(c_1+c_2\ln r\big),\\
    &r^{\frac{p-m}{2}}\big(c_1r^{\sqrt{\mu}}+c_2 r^{-\sqrt{\mu}}\big)
\end{cases}
\end{align*}
corresponding to the cases \(\mu<0,\mu=0,\mu>0\), where
\begin{align*}
 \mu=\frac{(m-p)^2}{4}-\ell(\ell+m-2)+\epsilon.   
\end{align*}
Hence, we can deduce that we have to require that \(\mu<0\) for some small \(\epsilon\) which is equivalent to the condition \eqref{eq:condition-p-harmonic-instability}.
We can conclude that
\begin{align*}
v(r)=
\begin{cases}
    r^\frac{p-m}{2}\sin(\sqrt{-\mu}\ln r), &\text{ for } r_0<r\leq 1,\\
    0, &\text { for } r\leq r_0
\end{cases}
\end{align*}
is the desired solution of \eqref{eq:ode-v-p-harmonic} completing the proof.
\end{proof}

As a next step towards completing the proof of Theorem \ref{thm:p-harmonic} we show, using a weighted Hardy inequality, that the generalized equator map \(u_\star^{(\ell)}\colon B^m\to\s^{m^{\ell}}\) is stable
assuming a certain inequality involving \(m,\ell,p\).

\begin{Lem}
For all \(\eta\in W_0^{1,2}(B^m,\R^{n+1})\)
and any \(p<m\) the following inequality holds
\begin{align}
\label{eq:hardy-weighted}
\int_{B^m}\frac{|\eta|^2}{r^p}dx\leq\frac{4}{(m-p)^2}
\int_{B^m}\frac{|d\eta|^2}{r^{p-2}} dx.    
\end{align}
\end{Lem}

\begin{proof}
For a straightforward proof we refer to the introduction of
\cite{hardy} and references therein.
See also \cite[Cor. 2.2]{MR1616905} for how to restrict
Hardy inequalities from Euclidean space to the unit ball.
Moreover, for radial functions the statement is established in
\cite[Lemma 1]{MR1743430}.

Again, these inequalities are all formulated for functions on 
\(B^m\)
but also hold for functions taking values in \(\R^{n+1}\)
following the arguments presented in \cite[p. 272]{MR4817500}.
\end{proof}

\begin{Prop}
Consider \(u_\star^{(\ell)}\colon B^m\to\s^{m^{\ell}}\) defined in \eqref{equator-generalized} as a $p$-harmonic map. Moreover, suppose that \(2\leq p<m\).
If the inequality
\begin{align}
\label{eq:condition-p-harmonic}
\frac{(m-p)^2}{4\ell(\ell+m-2)} \geq 1
\end{align}
is satisfied, then \(u_\star^{(\ell)}\) is a stable \(p\)-harmonic map.
\end{Prop}

\begin{proof}
We employ the second variation formula \eqref{eq:sv-p-energy} and insert the properties of the map \(u_\star^{(\ell)}\).
Then, from \eqref{eq:sv-p-energy} we get
\begin{align*}
\hess E_p(u^{(\ell)}_\star)(\eta,\eta)\geq
\ell^{\frac{p-2}{2}}(\ell+m-2)^\frac{p-2}{2}\int_{B^m}\frac{|d\eta|^2}{r^{p-2}}dx
-\ell^{\frac{p}{2}}(\ell+m-2)^\frac{p}{2}\int_{B^m}\frac{|\eta|^2}{r^{p}}dx,
\end{align*}
where we estimated the last term of \eqref{eq:sv-p-energy} by zero
as we are assuming that \(p\geq 2\).
Now, employing the weighted Hardy-inequality \eqref{eq:hardy-weighted} 
we find 
\begin{align*}
\hess E_p(u^{(\ell)}_\star)(\eta,\eta)\geq    
\ell^{\frac{p}{2}}(\ell+m-2)^\frac{p}{2}
\bigg(\frac{(m-p)^2}{4\ell(\ell+m-2)}-1\bigg)\int_{B^m}\frac{|\eta|^2}{r^{p-2}}dx,
\end{align*}
from which we get the result.
\end{proof}

Finally, we prove the last part of the statement of Theorem \ref{thm:p-harmonic} concerning the minimizing properties
of the generalized equator map \(u_\star^{(\ell)}\colon B^m\to\s^{m^\ell}\).

The next Proposition gives a refined statement on the stability
of the \(p\)-harmonic map \(u_\star^{(\ell)}\) by showing that it is not only a stable solution but even a minimizer of the \(p\)-energy \eqref{p-energy}.

\begin{Prop}
Consider \(u_\star^{(\ell)}\) defined in \eqref{equator-generalized} as a $p$-harmonic map. Moreover, suppose that \(2\leq p<m\).
If \eqref{eq:condition-p-harmonic} holds then \(u_\star^{(\ell)}\)
is a minimizing \(p\)-harmonic map.
\end{Prop}

\begin{proof}
Let \(w\in W^{1,2}(B^m,\s^{m^{\ell}})\) be arbitrary but with boundary value \(w=u_\star^{(\ell)}\)
on \(\partial B^m\). Using the weighted Hardy inequality \eqref{eq:hardy-weighted} we deduce
\begin{align*}
\int_{B^m}|\nabla u_\star^{(\ell)}|^{p-2}|\nabla(u_\star^{(\ell)}-w)|^2dx
\geq \int_{B^m}|\nabla u_\star^{\ell}|^p(u_\star^{(\ell)}-w)^2dx.
\end{align*}
By direct calculations we obtain the following two identities
\begin{align*}
\int_{B^m}|\nabla u_\star^{(\ell)}|^{p-2}|\nabla(u_\star^{(\ell)}-w)|^2dx
=&\int_{B^m}|\nabla u_\star^{(\ell)}|^pdx
-2\int_{B^m}|\nabla u_\star^{(\ell)}|^{p-2}\langle\nabla u_\star^{(\ell)},\nabla w\rangle dx \\
&+\int_{B^m}|\nabla u_\star^{(\ell)}|^{p-2}|\nabla w|^2dx,\\
\int_{B^m}|\nabla u_\star^{(\ell)}|^{p}(u_\star^{(\ell)}-w)^2dx
=&\int_{B^m}|\nabla u_\star^{(\ell)}|^{p}(2-2\langle u_\star^{(\ell)},w\rangle) dx.
\end{align*}
Since \(u_\star^{(\ell)}\) is a \(p\)-harmonic map we get
\begin{align*}
\int_{B^m}|\nabla u_\star^{(\ell)}|^{p-2}\langle\nabla u_\star^{(\ell)},\nabla\beta\rangle dx=
\int_{B^m}|\nabla u_\star^{(\ell)}|^{p}\langle u_\star^{(\ell)},\beta\rangle dx
\end{align*}
for all \(\beta\in W^{1,p}_0(B^m,\R^{m^{\ell}+1})\). Choosing \(\beta=u_\star^{(\ell)}-w\) this gives
\begin{align*}
\int_{B^m}|\nabla u_\star^{(\ell)}|^{p-2}\langle\nabla u_\star^{(\ell)},\nabla w\rangle dx=
\int_{B^m}|\nabla u_\star^{(\ell)}|^{p}\langle u_\star^{(\ell)},w\rangle dx.    
\end{align*}
Using the assumption \eqref{eq:condition-p-harmonic} and combining the previous inequalities
we get
\begin{align*}
\int_{B^m}|\nabla u_\star^{(\ell)}|^{p-2}|\nabla w|^2 dx\geq
\int_{B^m}|\nabla u_\star^{(\ell)}|^pdx.
\end{align*}
Moreover, employing Hölder's inequality we find
\begin{align*}
\int_{B^m}|\nabla u_\star^{(\ell)}|^{p-2}|\nabla w|^2\leq
\bigg(\int_{B^m}|\nabla u_\star^{(\ell)}|^{p}dx\bigg)^\frac{p-2}{p}
\bigg(\int_{B^m}|\nabla w|^{p}dx\bigg)^\frac{2}{p}.
\end{align*}
Combining the previous inequalities we thus arrive at
\begin{align*}
\int_{B^m}|\nabla u_\star^{(\ell)}|^pdx\leq \int_{B^m}|\nabla w|^pdx    
\end{align*}
for all \(w\in W^{1,p}(B^m,\s^{m^{\ell}})\) with boundary value \(w=u_\star^{(\ell)}\)
whenever \eqref{eq:condition-p-harmonic} holds.

From the Hardy inequality \eqref{eq:hardy-weighted} we deduce that the equality case is realized for \(w=u_\star^{(\ell)}\). Hence, under the assumption \eqref{eq:condition-p-harmonic} we have
\begin{align*}
\int_{B^m}|\nabla u_\star^{(\ell)}|^pdx<\int_{B^m}|\nabla w|^pdx     
\end{align*}
for all \(w\in W^{1,p}(B^m,\s^n)\) with boundary value \(w=u_\star^{(\ell)}\)
and \(w\neq u_\star^{(\ell)}\) in the interior.
Hence, we can conclude that \(u_\star^{(\ell)}\) is the unique minimizer of the \(p\)-energy \eqref{p-energy}
if \eqref{eq:condition-p-harmonic} holds, completing the proof.
\end{proof}

In the following remarks we connect the statement of
Theorem \ref{thm:p-harmonic} with the results that have already been established in the literature concerning the stability 
of the equator map.

\begin{Bem}
\begin{enumerate}
\item We can rewrite \eqref{eq:condition-p-harmonic} as follows
\begin{align}
\label{iq-m}
m>p+2\ell+2\sqrt{l}\sqrt{2\ell+p-2}.    
\end{align}
\item For \(\ell=1\) we reobtain \cite[Theorem A]{MR1743430}, 
where the classic Jäger-Kaul theorem for harmonic maps was extended to the case of \(p\)-harmonic maps.
More precisely, it was shown that the equator map \(u_\star^{(1)}\) is a minimizer of the \(p\)-energy \eqref{p-energy} if \(m\geq 2+p+2\sqrt{p}\), see also \cite{MR1647881}.
    \item For $\ell=1$,  Fardoun \cite{MR1626430} generalized 
    Theorem\,\ref{thm:p-harmonic}
    for maps taking values in ellipsoids.
\item 
In \cite{MR4579395} we established an existence result
for rotational \(p\)-harmonic self-maps of \(\s^m\).
Due to the rotational symmetry the \(p\)-harmonic map equation
reduces to an ODE which can be solved by a shooting method
under the assumption 
\begin{align*}
p<m<2+p+2\sqrt{p}    
\end{align*}
which matches \eqref{iq-m} in the case \(\ell=1\).
However, note that in \cite{MR4579395} the above condition was necessary to prove the existence of \(p\)-harmonic maps
whereas in Theorem \ref{thm:p-harmonic} it is central 
in the stability analysis of \(p\)-harmonic maps.
   
   \item Theorem\,\ref{thm:p-harmonic} can be restated in the following equivalent way:

\begin{Satz}
\label{thm:p-harmonic-2}
Assume that \(m>p\geq 2\). 
\begin{enumerate}
    \item If
        \begin{align*}
\ell>\frac{1}{2}(2-m+\sqrt{(m-2)^2+(m-p)^2}),
\end{align*}
     then the generalized equator map \(u_\star^{(\ell)}\colon B^m\to\s^{m^\ell}\) is an unstable \(p\)-harmonic map.
     \item If
        \begin{align*}
\ell<\frac{1}{2}(2-m+\sqrt{(m-2)^2+(m-p)^2}),
\end{align*}
     then the generalized equator map \(u_\star^{(\ell)}\colon B^m\to\s^{m^\ell}\) is the unique energy minimizer 
     of the \(p\)-energy \eqref{p-energy}.
\end{enumerate}
\end{Satz}

\item
In \cite[Remark $1$]{MR4757446} Nakauchi showed by long, explicit calculations that for
\begin{align*}
\ell\geq \frac{1}{2}(2-m+\sqrt{(m-2)^2+2(m-p)^2})
\end{align*}
the map \(u^{(\ell)}\) is unstable as a \(p\)-harmonic map.
Since 
\begin{align*}
\frac{1}{2}(2-m+\sqrt{(m-2)^2+2(m-p)^2})\geq\frac{1}{2}(2-m+\sqrt{(m-2)^2+(m-p)^2}),
\end{align*}
Theorem\,\ref{thm:p-harmonic} in particular improves the result of Nakauchi \cite{MR4757446}.    
   \item 
In \cite{MR1619840} it was proved that the equator map \(u_\star^{(1)}\colon B^m\to\s^{m}\) is \(p\)-energy minimizing for all \(p\in[m-1,m)\)
and the complementary case was addressed in \cite{MR1867937}, i.e. for \(m\geq 3\)
and \(1<p\leq m-1\) the map \(u_\star^{(1)}\colon B^m\to\s^{m}\) is the unique
minimizer of the \(p\)-energy \eqref{p-energy}.
\end{enumerate}
\end{Bem}

\appendix

\section{Proof of Lemma \texorpdfstring{\ref{app}}{a}}
In this appendix we prove Lemma\,\ref{app} by a brute force calculation. We expect that there exists a refined proof, but we were not able to establish the statement in a more elegant way.

\begin{proof}[Proof of Lemma\,\ref{app}]
Using the very definitions of $f_1$ and $f_2$ we get with the help of a computer algebra system that
\small
\begin{align*}
&f_1(\ell,s)-f_2(\ell,s)= -12006225 + 55583010 \ell + 644529087 \ell^2 - 2421377136 \ell^3 - 
 9700301278 \ell^4 \\& +12487858156 \ell^5 + 72563051494 \ell^6 + 
 94882488616 \ell^7 + 53978095723 \ell^8 + 13617646594 \ell^9 \\&+ 
 1475727899 \ell^{10} + 263739960 \ell^{11} + 52964100 \ell^{12} + 134390340 s + 
 1188352260 \ell s - 9358205616 \ell^2 s \\& -53535104072 \ell^3 s +
 82468999288 \ell^4 s + 670777476696 \ell^5 s + 1134928446416 \ell^6 s +  
 828852193984 \ell^7 s \\& +268985623908 \ell^8 s + 30497943172 \ell^9 s + 
 1148862144 \ell^{10} s + 457581480 \ell^{11} s + 1303670448 s^2 \\& - 
 10642384584 \ell s^2 - 120720353388 \ell^2 s^2 + 144787874944 \ell^3 s^2 + 
 2450530361488 \ell^4 s^2 + 5729642548976 \ell^5 s^2\\&  + 
 5589520852728 \ell^6 s^2 + 2486390322368 \ell^7 s^2 + 
 437882283680 \ell^8 s^2 + 15819743384 \ell^9 s^2 + 894836116 \ell^{10} s^2 \\& - 
 3642402384 s^3 - 123101112848 \ell s^3 - 34749970048 \ell^2 s^3 + 
 4173384480384 \ell^3 s^3 + 14832331134048 \ell^4 s^3\\&  + 
 20133076297696 \ell^5 s^3 + 12432504294400 \ell^6 s^3 + 
 3334091269376 \ell^7 s^3 + 291082689648 \ell^8 s^3\\&  + 3900335792 \ell^9 s^3 - 
 51483158224 s^4 - 255426616000 \ell s^4 + 3176624602464 \ell^2 s^4 + 
 20482984437760 \ell^3 s^4\\&  + 41075732638336 \ell^4 s^4 + 
 35735452614720 \ell^5 s^4  + 13867113842720 \ell^6 s^4 + 
 2075666251392 \ell^7 s^4\\&  + 77435615952 \ell^8 s^4 - 146586812800 s^5 + 
 731366779520 \ell s^5 + 14433759766528 \ell^2 s^5 + 
 47762510732544 \ell^3 s^5\\&  + 60993939304064 \ell^4 s^5 + 
 33813612386432 \ell^5 s^5 + 7594848719616 \ell^6 s^5 + 
 503719360512 \ell^7 s^5\\&  - 89987900416 s^6 + 4379984983808 \ell s^6 + 
 30133300439680 \ell^2 s^6 + 61633397659648 \ell^3 s^6 + 
 50136650047232 \ell^4 s^6\\&  + 16282468373760 \ell^5 s^6 + 
 1639823210112 \ell^6 s^6 + 359804487168 s^7 + 8997065145856 \ell s^7 + 
 35061994088448 \ell^2 s^7\\&  + 45067242383360 \ell^3 s^7 + 
 21502764900864 \ell^4 s^7 + 3144953498112 \ell^5 s^7 + 941458015488 s^8 + 
 9901678553600 \ell s^8 \\& + 23377417362176 \ell^2 s^8 + 
 17492946874368 \ell^3 s^8 + 3757798036992 \ell^4 s^8 + 1053226767360 s^9 \\& + 
 6202832655360 \ell s^9 + 8359807168512 \ell^2 s^9 + 
 2801413404672 \ell^3 s^9 + 640508858368 s^{10} + 2086916315136 \ell s^{10}\\&  + 
 1244059272192 \ell^2 s^{10} + 206616457216 s^{11} + 293232914432 \ell s^{11} + 
 27769409536 s^{12}.
\end{align*}
\normalsize
Now, we make repeated use of the fact that $s,\ell\geq 1$.
For example, we have $-12006225 + 55583010 \ell\geq 43576785\ell$. Proceeding this way, straightforward estimates yield
\begin{align*}
f_1(\ell,s)-f_2(\ell,s)\geq p(\ell,s),
\end{align*}
where $p$ is a non-zero polynomial in $\ell$ and $s$ with all coefficients non-negative. 
Thus $ p(\ell,s)>0$ and therefore
\begin{align*}
f_1(\ell,s)-f_2(\ell,s)>0,
\end{align*}
whence the claim.
\end{proof}

\bibliographystyle{plain}
\bibliography{mybib}

\end{document}